\documentclass[12pt]{amsart}
\begin{document}

\makeatletter
\let\over\@@over
\def\Box${ $\qed\ignorespaces}

\def\seq#1#2{#1_1,\ldots,#1_{#2}}
\def\pa#1#2{{\textstyle\partial #1\over\textstyle\partial #2}}
\def\pp#1{\mbox{$#1$-}}

\title {Polynomial Retracts and the Jacobian Conjecture} 

\subjclass{Primary  13B25, 13P10; ~Secondary 14E09.}

\author{Vladimir Shpilrain} 
\address{\hskip-\parindent Vladimir Shpilrain, 
Department of Mathematics, University of California, 
 Santa Barbara,  CA 93106, USA}
\email{shpil@math.ucsb.edu}

\author{Jie-Tai Yu} 
\address{\hskip-\parindent  Jie-Tai Yu,
Department of Mathematics, University of Hong Kong,  
Pokfulam Road, Hong Kong}
\email{yujt@hkusua.hku.hk}

\thanks{Research at MSRI is supported in part by NSF grant DMS-9022140.}

\begin{abstract}
Let $ K[x, y]$ be the polynomial algebra in two variables 
 over a field $K$ of characteristic $0$.  A subalgebra $R$ of  $K[x, y]$ is called a 
  retract if there is an idempotent homomorphism (a {\it retraction}, or {\it projection}) 
   $\varphi:   K[x, y] \to  K[x, y]$    such 
  that   $\varphi(K[x, y]) = R$. The presence of  other,
 equivalent, definitions 
of  retracts provides several different methods of studying them, 
  and brings together ideas from 
 combinatorial algebra, homological  algebra, and  algebraic 
geometry. In this paper, we 
   characterize all the  retracts of $ K[x, y]$ up to  an 
automorphism, and give  several 
 applications of this characterization, in particular,    to the
well-known Jacobian conjecture. 
  Notably, we prove that if  a  polynomial mapping $\varphi$  of $K[x,y]$ has 
  invertible Jacobian  matrix {\it  and }    fixes a non-constant 
polynomial, then  $\varphi$  is
   an automorphism. 
\end{abstract} 

\maketitle

\section{Introduction }
 
 Let $ K[x, y]$ be the polynomial algebra in two variables 
 over a field $K$ of characteristic $0$. A subalgebra $R$ of  $K[x, y]$ is called a 
{\it retract} if it satisfies any of the following equivalent conditions: 
\medskip

 \noindent {\bf (R1)} There is an idempotent homomorphism  (a {\it
retraction}, or {\it projection}) 

  \noindent     $\varphi:   K[x, y] \to  K[x, y]$    such 
  that   $\varphi(K[x, y]) = R.$ 
\smallskip

 \noindent {\bf (R2)} There is a homomorphism   $\varphi:   K[x, y]
\to R$ that fixes every element of 

 \noindent $R$. 
\smallskip

 \noindent {\bf (R3)}  $ K[x, y] = R \oplus I$  ~for some  ideal   
$I$ of the algebra   $K[x, y]$. 
\smallskip

 \noindent {\bf (R4)} $ K[x, y]$ is a projective extension of  $R$ in the category of  $ K$-algebras. 
   In other words, there is a split  exact sequence  $1  \to  I \to 
K[x, y] \to R \to 1$, 
 where  $I$ is the same ideal as in (R3) above. 

\medskip

 Examples: $ K$; $ K[x, y]$; any  subalgebra of the form $K[p]$, where $p \in  K[x, y]$ is a 
{\it coordinate}  polynomial (i.e.,  $K[p, q] = K[x, y]$ for some  polynomial  $q \in  K[x, y]$). 
 There are other, less obvious, examples of retracts: if $~p = x + x^2y$, 
then  $K[p]$ is  a retract of  $ K[x, y]$, but $p$ is not  coordinate since it has a fiber 
  ${\{}p = 0{\}}$ which is reducible,  and therefore is not isomorphic to a line. 
\medskip
  
 The very presence of several equivalent definitions 
 of  retracts shows how natural these objects are. Later on, we shall also comment on  a very natural 
geometric meaning of retracts. 
\smallskip

 In \cite{Costa}, Costa has proved that every proper retract of  $ K[x, y]$ (i.e., a one different 
from $K$   and  $ K[x, y]$)  has the form  $K[p]$  for some polynomial $p \in  K[x, y]$,  ~i.e., 
 is isomorphic to a polynomial $ K$-algebra in one variable. 
 A natural problem now is to  characterize somehow those polynomials $p \in  K[x, y]$ that 
 generate a retract of  $ K[x, y]$. Since the image of  a retract under  any automorphism of $ K[x, y]$ 
is again  a retract, it would be reasonable to characterize retracts up to  an automorphism of 
$ K[x, y]$, i.e., up to a ``change of  coordinates''. 
 We give an answer to this problem by proving the following 

\medskip

\noindent {\bf Theorem 1.1.} Let $K[p]$ be a retract of  $ K[x, y]$. 
There is 
 an automorphism $\psi$ of $ K[x, y]$ that takes the  polynomial 
$~p~$ to $~x + y \cdot q$ ~for some 
polynomial $q = q(x,y)$. A retraction for $K[\psi(p)]$ is given then
by $~x \to ~x + y \cdot q;   ~y \to 0$. 

\medskip 

 Our proof of this result is based on the well-known  Abhyankar-Moh
theorem \cite{AbMoh}. We   note in passing that  Theorem 1.1 also
yields a characterization of retracts of a free 
 associative algebra  $ K\langle x, y\rangle$ (see Theorem 2.1  in
the next Section 2) if one uses a 
   natural lifting. In Section 2, we also make an  
observation on  retracts of a polynomial algebra in arbitrarily many
variables 
  (Proposition 2.2). 

\smallskip 
 Theorem 1.1 yields  another useful 
  characterization of retracts of $ K[x, y]$: 
\medskip

\noindent {\bf Corollary 1.2.} A  polynomial $~p \in  K[x, y]$ 
generates a retract of  $ K[x, y]$ 
  if and only if there is  a  polynomial mapping of $K[x,y]$ that 
takes  $~p~$ to $~x$.  The ``if" part is
actually valid for a  polynomial algebra  in arbitrarily many
 variables. 

\medskip

 Although the form to which any retract can be reduced by Theorem 
1.1  might seem rather general,  it is in fact quite restrictive,
 and has several interesting applications, in particular, 
  to the notorious 
\medskip 

\noindent {\bf Jacobian conjecture.}\cite{Keller}
~If for a pair of polynomials  $p,q \in  K[x, y]$, the 
 corresponding Jacobian matrix is invertible, 
 then  $ K[p, q] =  K[x, y]$. 

\medskip
 For a survey and background on this problem, the reader is 
referred to ~\cite{BCW}.  
\smallskip 

 To explain how  retracts of  $ K[x, y]$ and the Jacobian conjecture 
are connected, we need one  more 
\medskip

\noindent {\bf Definition.}  A  polynomial $p \in  K[x, y]$  has a  
{\it unimodular   gradient } if there are   
   $u,v \in  K[x, y]$   such 
  that   $ \pa{p}{x} \cdot u +   \pa{p}{y} \cdot v = 1$. 

\medskip

 If the   ground field $~K$ is 
algebraically closed, then this  is equivalent, by Hilbert's
Nullstellensatz, to the  gradient being nowhere-vanishing. 
\smallskip 

 Now we are ready to establish a link between  retracts of   $ K[x,
y]$  and the Jacobian conjecture  by means of the following 
\medskip

\noindent {\bf Conjecture ``R''.}  If a  polynomial $p \in  K[x, y]$ 
has a  unimodular gradient,  then  $K[p]$ is  a retract of  $ K[x,
y]$.

\medskip
 
\noindent {\bf   Theorem 1.3.}    Conjecture ``R'' implies the 
Jacobian conjecture. 

\medskip

 In fact,  Conjecture ``R'' is stronger than the Jacobian conjecture since there are  polynomials 
that have a  unimodular gradient, but do not  have a Jacobian mate. 
The simplest example is  $~p = x + x^2y$; ~this polynomial  has a 
unimodular gradient, but this gradient cannot be  completed to an
invertible $2\times 2$  matrix by the gradient of any  other
polynomial.  
\smallskip 

 Therefore, a statement {\it equivalent} to the Jacobian conjecture 
would be: if $(p,q)$ is  a Jacobian pair (i.e., the 
 corresponding Jacobian matrix is invertible), then  $K[p]$   is  a 
retract of  $ K[x, y]$.
 This statement       is formally weaker than the Jacobian
conjecture, but is actually equivalent to it by our  Theorem 1.3. 
\medskip

 There are several different 
ways of proving  Theorem 1.3  using the characterization of retracts 
given in  Theorem 1.1; we believe that the proof we give here (using
Newton polygons) is  particularly  elegant. 
\smallskip 

 As a corollary, we show that the Jacobian conjecture is 
equivalent to   other (formally) much weaker statements:  
\medskip

\noindent {\bf Corollary 1.4.} Either of the 
 following claims is 
equivalent to the Jacobian conjecture. Suppose $(p,q)$ is 
 a Jacobian pair. Then: 
\smallskip 

\noindent {\bf (i)} for some  coordinate polynomial $g$, 
$~K[x,y] = K[p] ~+ < g >$,  
~where $~< g >$ ~is the ideal of $ K[x,y]$ generated by 
 the polynomial $g$.
\smallskip 

\noindent {\bf (ii)} for some  coordinate polynomial $g$, 
 $~K[x,y] = K[g] ~+ < p >$. 

\medskip

Corollary 3.1 in Section 3 also seems interesting since it improves
several known partial results on two-variable 
 Jacobian conjecture.

\medskip 

  Conjecture ``R''  has the following 
 geometric interpretation. Suppose 
we have a polynomial  $p \in  \mathbf{C}[x, y]$ with  
nowhere-vanishing gradient. Then every fiber  ${\{}p = c, ~c \in
\mathbf{C}{\}}$  is a non-singular curve on the complex plane. Fix a
particular fiber, say, ${\{}p = c_0{\}}$. 
 If there is another polynomial $q \in  \mathbf{C}[x, y]$ such 
  that every fiber of $q$ intersects the curve ${\{}p = c_0{\}}$ 
 transversally (i.e., without  common tangent line) at exactly one
point, then we can arrange a  geometric  projection of the 
  plane onto the curve ${\{}p = c_0{\}}$ by sliding a point on a 
fiber  ${\{}q = c{\}}$ 
 toward the intersection  point of  ${\{}q = c{\}}$ with ${\{}p =c_ 0{\}}$. This  geometric  projection 
will also be algebraic, i.e., $\mathbf{C}[p]$  will be a   retract 
of  $\mathbf{C}[x, y]$ in this case. 

\smallskip 

 As far as the Jacobian conjecture is concerned, the Jacobian
condition ensures transversal  intersection 
 of any fiber of $p$ with any fiber of $q$ if  $(p,q)$ is 
 a Jacobian pair.  Moreover, by a result of Kaliman \cite{Kaliman}, 
we can restrict our attention to the situation  where all the fibers
of $~p$  and $~q$ are irreducible, i.e., are connected curves. The
only problem  is to show that there is a fiber of $~p$ which 
intersects every  fiber of $~q$ at {\it  exactly} one point, 
   i.e., that for some  $c_0  \in \mathbf{C}$, the system of
equations   $p(x,y) = c_0; ~q(x,y) = c$  has a unique solution for
any $~c \in \mathbf{C}$. Thus, we have: 

\medskip

\noindent {\bf Corollary 1.5.} (cf. \cite{Gwos}). Suppose $\varphi$ 
is  a  polynomial mapping of $\mathbf{C}[x,y]$ 
 with invertible Jacobian  matrix. If $\varphi$ is injective on 
some  line, then $\varphi$ is 
  an automorphism. 
\medskip

 We also note here 
          that the aforementioned 
 result of Kaliman \cite{Kaliman} calls for an example 
of a non-coordinate polynomial $~p \in \mathbf{C}[x,y]$ that has a
non-vanishing 
 gradient 
  and all of whose fibers are irreducible. If there is no 
polynomial like that, then the two-variable Jacobian conjecture 
is true. 
\medskip 

 Another application of retracts to the Jacobian conjecture 
(somewhat indirect though) is based on  the
``$\varphi^{\infty}$-trick'' familiar in combinatorial group theory
(see  \cite{Turner}).  For a polynomial mapping  $\varphi : K[x,y]
\to K[x,y]$, ~denote by $\varphi^{\infty}(K[x,y]) = 
 \bigcap_{k=1}^{\infty}  \varphi^k(K[x,y])$ ~the {\it stable image} 
of $\varphi$. Then we have: 
\medskip

\noindent {\bf Theorem 1.6.}      Let $\varphi$ be a  polynomial 
mapping of $K[x,y]$. If the Jacobian  matrix of $\varphi$ is
invertible,  then either   $\varphi$  is   an automorphism, or 
$~\varphi^{\infty}(K[x,y]) = K$.  
\smallskip 

  Our proof of Theorem  1.6 
 is based on  recent results of Formanek \cite{Formanek} and 
 Connell-Zweibel \cite{Connell}. 
   
\smallskip 

 Obviously, if  $\varphi$ fixes  a  polynomial $p \in  K[x, y]$, 
~then $~p \in \varphi^{\infty}(K[x,y])$.
 Therefore, we have: 
\medskip

\noindent {\bf Corollary 1.7.} Suppose $\varphi$ is  a  polynomial 
mapping of $K[x,y]$ 
 with invertible Jacobian  matrix. If $\varphi(p) = p$  ~for some
non-constant polynomial 
$p \in  K[x, y]$, ~then  $\varphi$  is   an automorphism. 
\medskip

  This yields the following   re-formulation of the Jacobian
conjecture: if  $\varphi$ is  a  polynomial mapping of $K[x,y]$ 
 with invertible Jacobian  matrix, then for some automorphism 
$~\alpha$, the  mapping  $~\alpha  \cdot \varphi~$ fixes a
non-constant polynomial. \\

\section{Retracts of $\mathbf{~K[x, y]}$}

 We start with 
\smallskip 

\noindent {\bf Proof of Theorem 1.1.}  Let $K[p]$, $~p \in  K[x, y]$, 
 be a retract of $ K[x, y]$ ~(note that by a result of Costa
\cite{Costa}, every retract of $ K[x, y]$ has this form). Without
loss of generality, we may assume that $~p$ has zero constant term.

 Let the corresponding retraction be given by  $\phi: x \to
q_1(p); ~y \to q_2(p)$ ~for some one-variable polynomials $~q_1,
q_2$. Again, we may assume that both $~q_1$ and $~q_2$ have 
zero constant term. 

 Then, since $\phi$ is a retraction,   polynomials $~q_1(p),
q_2(p)~$  should generate $K[p]$. Suppose $~q_1$ and $~q_2$ have 
 degree $n \ge 1$ and $m \ge 1$, respectively. Then, by 
Abhyankar-Moh theorem \cite{AbMoh}, either $~n$ divides $~m$, or 
 $~m$ divides $~n$. Suppose $~deg(q_1) = k \cdot deg(q_2)$ ~for
some integer $~k \ge 1$. 

 Then we make the following change of coordinates: $x \to 
 \widetilde x  = x - c \cdot  y^k; ~y \to \widetilde y = y$, where 
the coefficient  $c \in K^\ast$ is chosen so that in the polynomial 
 $~q_1 - c \cdot q_2^k$, the leading terms would cancel out. 
\smallskip 

  In these new coordinates, our retraction acts as follows: 
 $\phi: \widetilde x \to q_1 - c \cdot  q_2^k = \widetilde q_1$; 
 $~\widetilde y \to q_2 = \widetilde q_2$. Polynomials 
$~\widetilde q_1$ and $~\widetilde q_2$  are easily seen to be
   another generating set of $~K[p]$, but the sum of degrees of 
$~\widetilde q_1$ and $~\widetilde q_2~$ is less than that of 
$~q_1~$ and $~q_2$. 
 
 Continuing this process, we shall eventually arrive at a pair of 
polynomials one of which is zero. Denote the other one by $~h$;
then we must have $K[h] = K[p]$, i.e., $~h = c \cdot p$  for some 
$c \in  K^\ast$. 

 Thus, we have shown that for some automorphism $~\psi \in  Aut(K[x,
y])$  (``change  of  coordinates"),  the composition $~\phi \psi$ 
~takes  $~x~$ to   $~p$,  and  $~y~$ to 0. It follows that 
 in these coordinates, 
 $~p(x,y) = x + y \cdot q(x,y)$ ~for some polynomial $~q$. This
completes the proof. 
$\Box$ 
\medskip 
 
\noindent {\bf Proof of Corollary 1.2.} (1) Suppose $~p \in  K[x, y]$ 
generates a retract of  $ K[x, y]$. Then, by Theorem 1.1, for some 
 automorphism $~\psi \in  Aut(K[x,y])$, the   polynomial $~\psi(p)$ 
has the form $~x + y \cdot q(x,y)$. Let $~\phi $ be a mapping of 
 $ K[x, y]$ that takes $~x $ to $~x$; $~y$  to 0. Then
$~\phi(\psi(p)) = x$. 
\smallskip 

\noindent (2) We are going to prove the ``if" part for a  polynomial 
algebra $ K[x_1,...,x_n]$ in arbitrarily many variables.

 Let $~\varphi(p) = x_1$.  Consider  the   following  mapping of
$K[x_1,...,x_n]$: 
  $~\psi: x_1 \to p; ~x_i \to 0$, $~i = 1,...,n$. Then  
 \begin{eqnarray}
 \psi(\varphi(p)) =p. 
\end{eqnarray} 
   Denote $~\varrho = \psi\varphi$. Then, by (1),
$~\varrho(p) = p$, which means  $~\varrho$   fixes every element of
$K[p]$.    Also, it is clear that  $~\varrho(K[x_1,...,x_n]) = K[p]$. 
 Therefore, $~\varrho$  is a retraction of $ K[x_1,...,x_n]$, and 
 $K[p]$ a retract. $\Box$ 
\medskip 

 Now we  give a characterization of retracts 
of a free  associative algebra  $ K\langle x, y\rangle$: 
\medskip 

\noindent {\bf Theorem 2.1.}   Let $R$ be a proper retract of   
$K\langle x, y\rangle$. There is 
 an automorphism $\psi$ of $K\langle x, y\rangle$ that takes $R$ 
        to $K\langle v\rangle = K[v]$ for some element  $~v$ of the
 form $~x + y \cdot q(x,y) + w(x,y)$, ~where  $~w(x,y) $ belongs to
the commutator ideal of $K\langle x, y\rangle$. 
\medskip 
 
\noindent {\bf Proof.} First of all,   every element 
of the given form generates a retract of   $K\langle x, y\rangle$. 
 Indeed, the corresponding retraction is 
 given by  $~x \to v; ~y \to 0$ ~(this mapping will take $w$ to 0 
as well).
\smallskip 

  Now we are going to show that every retract of   
$K\langle x, y\rangle$ has the form $~K\langle v\rangle = K[v]$ 
for some element  $~v \in K\langle x, y\rangle$. From the definition 
(R1) of a retract, we see that every retract of $K\langle x, y\rangle$ 
  can be generated by two elements ($\varphi(x)$  and $\varphi(y)$).

 Another easy observation is that if $R$ is a retract of 
  $K\langle x, y\rangle$ with the corresponding retraction $\phi$, 
then $R^{\alpha}$  is a retract of $ K[x, y]$, and the corresponding 
retraction is $\phi^{\alpha}$. Here  $R^{\alpha}$ denotes the image
of $R$ under the natural abelianization mapping  $~\alpha: K\langle
x, y\rangle \longrightarrow K[x, y]$ ~(note that the kernel of this 
mapping  is the commutator ideal of $K\langle x, y\rangle$),  and 
 $\phi^{\alpha}$ is the natural abelianization of $\phi$. 
 The best way to see it is to apply the  definition (R2) of a
retract. 
\smallskip 

 Upon combining these two observations with what we know about 
 retracts  of  $ K[x, y]$,  we see that generators of 
our retract $~R$ must be of the   form  $~v_1(x,y) = q_1(x,y) + 
w_1(x,y)$ ~and $~v_2(x,y) = q_2(x,y) + w_2(x,y)$, ~where  $w_1,
w_2$ belong to the commutator ideal of $K\langle x, y\rangle$, 
 ~and $q_1^{\alpha} = h_1(p), ~q_2^{\alpha} = h_2(p)$, 
 ~where  $p = p(x,y) \in  K[x, y]$ is  a polynomial   
  that   generates a retract of  $ K[x, y]$.
\smallskip 

 Moreover, since $R^{\alpha}$ should be equal to $K[p]$, we 
should have  $ K[h_1(p),h_2(p)] = K[p]$, so that we are in a
position to apply Abhyankar-Moh theorem. Repeating the argument from
the proof of Theorem 1.1, we see that after applying an appropriate 
 automorphism $~\psi \in  Aut(K[x, y])$  (``change  of 
coordinates"), our pair $(h_1(p),h_2(p)$ becomes  
$(\psi(h_1(p)),0)$. 
\smallskip 

 Since every automorphism of $ K[x, y]$ can be lifted to an  
automorphism of $K\langle x, y\rangle$ (see  \cite{Cohn}, Theorem
8.5), we have shown that, after applying an automorphism if
necessary, our retract $~R$  can be generated by two elements, 
one of which (call it $~u$) 
belongs to the commutator ideal of  $K\langle x,
y\rangle$, and the    other one (call it $~v$) is of the   form  $~x
+  y \cdot q(x,y) + w(x,y)$, ~where  $~w(x,y)$ belongs to the
commutator ideal of $K\langle x, y\rangle$,  and the corresponding 
retraction is  $~\phi: x \to v; ~y \to u$.

  Then, since $\phi^{\alpha}$ annihilates  $~y^{\alpha}$ (see
above),  we should have  $~u = \phi(y)$ belong to the commutator
   ideal of  $K\langle x,y\rangle$. 

 On the other hand, $~\phi$ should fix every element of $~R$, in 
particular, the element $~u$.  Suppose $~u \ne 0$. 
Write $~u~$ in the form $~u = 
    m_1 + \widetilde u$, ~where    $~m_1$ is the sum of  lowest
  degree terms, ~i.e., any monomial of $~\widetilde u$ 
 has degree greater than that of monomials in $~m_1$. Then, since 
 $~u~$ belongs to the commutator 
   ideal of  $K\langle x,y\rangle$, every monomial in $~m_1$ 
 depends on $~y$. Therefore, the image of any monomial in $~m_1$ 
under the endomorphism  $~\phi$ is a sum of monomials whose
 degree is greater than that of monomials in $~m_1$ (this sum
might be equal to zero, but it does not affect the argument). 
 The same applies to   monomials of $~\widetilde u$.  
\smallskip 

 Thus, there is no way we can have $~\phi(u) = u$; this
contradiction shows that  $~u = 0$, which completes the  proof of 
 the theorem. $\Box$ 
\medskip 

We conclude this section with an   observation on 
retracts of a polynomial algebra in arbitrarily many variables.  

\medskip 

\noindent {\bf  Proposition 2.2.} Let $~R$ be a proper retract of 
 $K[x_1,...,x_n]$ generated by polynomials $~p_1, ..., p_n$, 
 $~n \ge 2$. Then  $~p_1, ..., p_n$ are algebraically dependent.
\medskip 
 
\noindent {\bf Proof.} Let $~\phi: K[x_1,...,x_n] \to R$ ~be a 
retraction, so that $~\phi(R) = R$. In particular, $~\phi$ 
restricted to $~R~$  is an automorphism of $~R$. If 
polynomials $~p_1, ..., p_n$ ~were algebraically independent, then  
 a result of \cite{Connell} would imply that $~\phi$ is an 
automorphism of  $K[x_1,...,x_n]$ as well. In that case, we would 
have $~R = K[x_1,...,x_n]$, which means $~R~$ is not a proper
retract, hence a contradiction. $\Box$ 
\medskip 

 We note that it is an open problem whether or not any retract of 
 $K[x_1,...,x_n]$, $~n \ge 3$, ~can be generated by 
 algebraically independent polynomials. If this is   the case, it 
would imply the positive solution of the
well-known {\it cancellation problem} -- see \cite{Costa} for 
      discussion. \\

\section {The Jacobian conjecture}

 First we prove that our Conjecture ``R'' implies the 
Jacobian conjecture. 
\medskip 
 
\noindent {\bf Proof of Theorem 1.3.}  Let $\varphi: x \to
p(x,y); ~y \to q(x,y)$ ~be a 
polynomial mapping of $ K[x, y]$ with invertible Jacobian  matrix. 
Then, in particular, the polynomial $p(x,y)$ has a   unimodular  
gradient. 
 \smallskip 

 If we assume that our  Conjecture ``R'' is true, then $p(x,y)$ 
generates a retract of  $ K[x,y]$, hence  by  Theorem 1.1, for some 
 automorphism $~\psi \in  Aut(K[x,y])$,   polynomial $~\psi(p)$ 
has the form $~x + y \cdot h(x,y)$.  Therefore, upon combining 
 $\varphi$ with $~\psi$ if necessary, we may assume that $p(x,y)$ 
itself has this form. 
\smallskip 

 Now we appeal to a result of \cite{Lang1}    concerning 
 Newton polygons of a Jacobian pair  $(p(x,y), q(x,y))$.  
 The Newton polygon of a polynomial $~f = f(x,y) = \sum a_{ij} x^i
y^j$  ~is the convex hull of $~\{(i,j) \vert $ $a_{ij}
\neq 0 \}\cup \{(0,0) \}$. 
 \smallskip 

 The result of \cite{Lang1} we need is that if $(p(x,y), q(x,y))$ 
 is a Jacobian pair, but $~\varphi: x \to p(x,y); ~y \to q(x,y)$ is 
not an automorphism, then Newton polygons of 
$p(x,y)$ and $q(x,y)$  are radially similar. 
\smallskip 

 Now, by way of contradiction, suppose $~\varphi$ is 
not an automorphism. Look at the Newton polygon of $p(x,y)$. We see
       that it has an edge of length 1, namely, the one between 
the vertices (0,0)  and (1,0). It follows that the similarity ratio 
for Newton polygons of $p(x,y)$ and $q(x,y)$ is an integer; in 
particular, $q(x,y)$ has the form $~\sum_{i=1}^{k} c_i x^i + y\cdot
f(x,y)$ for some $c_i \in K$, $~c_k \ne 0$. 
 \smallskip 

 Then replace the pair $(p, q)$ with  
 $~(p, ~q - c_k p^k)$. This
new pair clearly has the same   properties as $(p, q)$ does: it is a
Jacobian pair, but the  corresponding mapping is not an
automorphism. However, the highest
 degree of monomials of the form $~x^m$ in the second polynomial 
has been decreased. 
 \smallskip 

 Therefore, we can 
repeatedly apply our argument (note that $p(x,y)$ does not change), 
until we get a pair $~(p, g)$, ~where  $~g$  has no monomials of the 
form $~x^m$. But in that case, the Newton polygon of $~g$  has no 
 edges along the $x$-axis, hence Newton polygons of $~p$ and
$~g$ cannot be radially similar. This contradiction completes the 
proof of  Theorem 1.3. $\Box$ 
\medskip

\noindent {\bf Proof of Corollary 1.4. (i)} Without loss of
generality, we may assume that both $~p$ and $~g$ have 
zero constant term. Moreover, upon applying an automorphism to 
all polynomials under consideration if necessary, we may assume 
 $~g = x$.
 \smallskip 

First we show that if 
$ K[x,y] = K[p] ~+ < x >$, ~then the sum is 
 actually direct. 
 \smallskip 

 By way of contradiction, suppose we have 
 \begin{eqnarray}
 x \cdot u = \sum_{i=1}^m c_i  \cdot p^i 
\end{eqnarray} 

\noindent  for some non-zero 
polynomial $~u = u(x,y)$ and constants
$~c_i  \in K$. Since the left-hand side of (2) is  divisible by
$~x$, the right-hand side  should be divisible by $~x$, too. 
 This is only possible if $~p$ itself is  divisible by $~x$, 
 but in that case, the Newton polygon of $p$ would not have an edge 
 along the $~y$-axis, which contradicts $p$ having a Jacobian mate 
 (see e.g. \cite{Lang1}). 

 \smallskip 

 Thus, we have shown that $ K[x,y] = K[p] \oplus \langle x
\rangle$. By the definition (R3) of a retract, this  implies 
 $ K[p]$  is a retract of $ K[x,y]$. Therefore, by our 
  Theorem 1.3, $ K[p,q] = K[x,y]$. 
$\Box$ 
\medskip 

\noindent {\bf (ii)} Arguing as in (i), we get $~p \cdot u =
\sum_{i=1}^m c_i  \cdot x^i$.   
 In this case, $ ~p$ cannot depend on $~y$, so $ ~p = p(x)$; but
then  the equality $ K[x,y] = K[x] + \langle p \rangle$ is not 
possible.

 Thus, $ K[x,y] = K[x] \oplus \langle p \rangle$. This clearly 
  implies $ ~p = c \cdot y$ ~for some $c \in K^\ast$. The result
follows. $\Box$ 
 \medskip 

\noindent {\bf  Corollary 3.1.} Let $K$ be an algebraically closed 
 field, and let $~(p, q), ~p,q \in K[x,y]$, be a Jacobian pair. 
 Suppose  $~p$ 
 has the form  $~p = x +  g$ ~for some 
 $~g \in K[x,y]$, which is  divisible by a homogeneous 
polynomial.  Then  $K[p,q] = K[x,y]$. 
\medskip 
 
\noindent {\bf Proof.} If $~p$ is linear, then we are done. 
Suppose $~p$ is  non-linear,  and suppose $~g$ is  divisible by a 
 non-constant homogeneous polynomial $~h$. 

If $~h$ is  divisible by $~x$, then $~p$ itself is  divisible by 
$~x$. In that case, the Newton polygon of $p$ does 
 not have an edge 
 along the $~y$-axis, which contradicts $~p$ having a Jacobian mate 
 (see  \cite{Lang1}). 

 If $~h$ is  divisible by $~y$, then $K[p]$ is a retract of        
 $K[x, y]$. If not, then there is $~c \in K^\ast$  such 
  that  the homomorphism $~x \to p; ~y \to cp$  takes $~h$ to 0. 
 This is obviously a retraction, so again, $K[p]$ is a retract of  
 $K[x, y]$. 

 Applying our Theorem 1.3 yields the result. $\Box$ 
 \medskip 

  Van den Essen and   Tutaj \cite{Essen} have shown
that if a Jacobian pair $~(p, q)$ is of the form $~(x + h_1, ~y +
h_2)$, ~where $both$   $~h_1$ and $~h_2$ are homogeneous 
polynomials, then $K[p,q] = K[x,y]$. It is notable that our 
 Corollary 3.1 not only relaxes the condition on the form of 
polynomials, but, most importantly, our condition is imposed on 
 $one$ polynomial only. 

\medskip 

 To prove Theorem 1.5, we need the following result of independent
interest: 
 \medskip 

\noindent {\bf Lemma 3.2.}  Let $\varphi: x \to
p(x,y); ~y \to q(x,y)$ ~be a 
polynomial mapping of $ K[x, y]$ with invertible Jacobian  matrix. 
 Suppose  $~\varphi(K[x, y])$  ~contains a coordinate polynomial.
 Then  $~\varphi$ is an 
 automorphism.  
\medskip 
 
\noindent {\bf Proof.} Upon composing $~\varphi$ with an 
 automorphism if necessary, we may assume that $~x \in  \varphi(K[x,
y])$. This implies $ K[p,q,x] = K[p,q]$, and, therefore, 
 $ K(p,q,x) = K(p,q)$, ~where  $ K(p,q)$ is the quotient field of 
$K[p,q]$.  

 On the other hand, by  a result of \cite{Formanek}, $(p, q)$ 
 being a Jacobian pair  implies $ K(p,q,x) = K(x,y)$. 

 Therefore, we have  $ K(p,q) = K(x,y)$, ~which by Keller's theorem
\cite{Keller} implies  $ K[p,q] = K[x,y]$.  $\Box$ 
\medskip 
 
\noindent {\bf Proof of Theorem 1.6.} Suppose $~\varphi$ is 
 \textsl{ not} an  automorphism.  Then, by Lemma 3.2, $~\varphi(K[x,
y])$   ~contains no coordinate polynomials. 
 In particular, the   degree of any polynomial in $~\varphi(K[x,
y])$  is at least 2. By using inductive argument, we are going to
show now that the   degree of any polynomial in $~\varphi^k(K[x,
y])$  is at least $~(k+1)$; ~this will imply $~\varphi^{\infty}(K[x,y]) 
= K$. 
 \smallskip 

 Consider an algebra  $~\varphi^k(K[x,y])$ for some $~k \ge 1$. 
 Since $~\varphi$ is injective (this is ensured by the Jacobian 
condition), polynomials $~\varphi^k(x)$  and $~\varphi^k(y)$ are 
 algebraically independent. 
 \smallskip 

 If $~\varphi$ restricted to $~\varphi^k(K[x,y])$ 
were  an 
 automorphism of   $~\varphi^k(K[x,y])$, then, by a result of 
 \cite{Connell}, $~\varphi$ would be an 
 automorphism of   $K[x,y]$, contrary to our assumption. (Note 
 that $~\varphi^k(K[x,y])$ is invariant under $~\varphi$ since, by 
 induction, $~\varphi^k(K[x,y]) \subseteq \varphi^{k-1}(K[x,y])$ 
 implies $~\varphi^{k+1}(K[x,y]) \subseteq \varphi^k(K[x,y])$). 
 Hence, $~\varphi \mid _{\varphi^k(K[x,y])}$ is 
 not an  automorphism of $~\varphi^k(K[x,y])$, and our previous 
 argument yields that $~\varphi^{k+1}(K[x,y])$ ~contains no 
coordinate polynomials of $~\varphi^k(K[x,y])$. 
 In particular, the   degree of any polynomial in 
$~\varphi^{k+1}(K[x,y])$ is greater than that of any coordinate 
polynomial in 
 $~\varphi^k(K[x,y])$. This completes the induction. $\Box$ 
\medskip

 A. van den Essen has pointed out to us that in the case $K =
\mathbf C$, a more complicated (geometric) proof of the result of  
  Theorem 1.6 was given by  Kraft \cite{Kraft}, who  also proved 
that if a polynomial 
mapping  $~\varphi$ of $~\mathbf C[x,y]$ is $not$ birational (i.e., 
 $~\varphi$ does not induce an  automorphism of the quotient  
field  $~\mathbf C(x,y)$), then 
 the   stable image  of $~\varphi$ is a retract of $\mathbf
C[x,y]$. 
 This yields a   natural question -- 
 is   the same true for $any$ polynomial 
mapping?  ~(cf. \cite{Turner}). \\

\section*{Acknowledgements}

 We are indebted to A.Campbell and D.Markushevich for numerous 
 insightful discussions. 
  The first author is grateful to Department of Mathematics of the 
University of Hong Kong for its warm hospitality during his visit
when this work was 
 initiated. A crucial part of this work 
was 
 done when the first author was visiting the Mathematical Sciences 
Research Institute in   Berkeley, and the second author 
was visiting the  Max Planck Institut f\"{u}r Mathematik in Bonn. 
The  warm hospitality and stimulating atmosphere of these
institutions  are greatly appreciated.

\medskip 

\end{document}